# COMING TO SEE FRACTIONS ON THE NUMBERLINE

Elisabetta Robotti*, Samuele Antonini**, Anna Baccaglini-Frank***

*Università della Valle d'Aosta, **Università di Pavia,

***Università di Modena e Reggio Emilia

*The aim of this paper is to present a didactical sequence that fosters the development of meanings related to fractions, conceived as numbers that can be placed on the number line. The sequence was carried out in various elementary school classes, containing students certified to have mathematical learning disabilities (MLD). Thus, our didactical aim is to make accessible meaning related to fractions to all students of the class, including MLD students, by means of a common didactical sequence, that is, a sequence proposed to the entire class. Our research is based on a range of different perspectives, from mathematics education to neuroscience and cognitive psychology. We discuss how such perspectives can be combined and provide the theoretical bases to design an outline of the didactical sequence, which allows us to implement and strengthen inclusive education.*

## INTRODUCTION AND LITERATURE

The concept of *fraction* is a very difficult one to master: frequently students are unable to reach an appropriate understanding of it, as described for example by Fandiño Pinilla (2007), and they can even come to fear fractions (Pantziara & Philippou, 2011). When children encounter fractions – typically in Italian school this happens in third grade (at 8-9 years of age) – it is the first time they have to treat sets of digits differently than those related to each other by decimal positional notation in positive integers. The numerator and denominator of a fraction are two numbers, each of which is bound by the rules that apply to positive integers, but that together represent a new, *single*, number. Learning to see the numerator and the denominator of a fraction together, as a single number is one of the most difficult – if not *the* most difficult – cognitive aspect of fractions (Bobis, Mulligan, & Lowrie, 2013).

Others disciplines, besides mathematics education, such as cognitive psychology and neuroscience, have also been very active in investigating the phenomena of (difficulties in) understanding of mathematics (incuded fractions), even if the different interested fields of research have not yet reached sufficiently common grounds for conducting scientific and interdisciplinary studies. In this paper, we consider some results from research in neuroscience and cognitive psychology to ground important design decisions taken during the elaboration of a teaching experiment constructed around the learning of fractions in primary school. In the following paragraphs we will illustrate reasons why it is important to learn fractions, both from a didactical (math education) point of view and from the perspective of cognitive science.

# IMPORTANCE OF FRACTIONS AND DIFFICULTIES ENCOUNTERED IN THEIR LEARNING

From the logical and epistemological points of view, the notion of *fraction* can be seen in different ways: as a linguistic representation of the decimal number obtained from the division indicated (but not calculated) by the number corresponding to the numerator and the one corresponding to the denominator; as an operator where the denominator indicates in how many equal parts a given unit is divided into (*each part is called a unit fraction*) and the numerator indicates the number of these to consider.

$$\frac{3}{4}u = 3 \; times \; \frac{1}{4}u \rightarrow 3 \times \frac{1}{4}u$$

Frequently, at least in Italian education, the conception of fraction as an operator is not explicitly identified as a rational number. Only when it is transformed into a decimal number is it placed on the number line.

From the point of view of learning mathematics, fractions constitute an important leap within domain of arithmetic because they represent a first approach towards the idea of extension of the set of Natural Numbers. In this sense, fractions need to assume a specific position on the number line (Bobis et al., 2013; Bartolini Bussi et al., 2013). Teaching the notion of fraction is, therefore, a quite delicate issue and it is ever so important to explore insightful ways of structuring didactical activities around it. In this respect, particularly insightful approaches have been provided, for example, by Bobis, Mulligan and Lowrie (2013). Even if certain basic aspects of the concept of fraction, particularly when seen as the perception of the variation of a ratio, seem to be innate (McCrink & Wynn, 2007), the learning of fractions presents obstacles, not only of a didactical nature. In fact, research in mathematics education (e.g., Bartolini Bussi et al., 2013), has shown how learning about not only the semantic aspects but also the lexical and syntactical ones of fractions involves the overcoming of different epistemological and cognitive obstacles such as:

- Assuming that the properties of ordering natural numbers can be extended to ordering fractions (e.g. assuming that the product/quotient of two fractions makes a greater/smaller fraction).

- Positioning fractions on the number line using the pattern of whole numbers (Iuculano & Butterworth, 2011).

From a cognitive point of view, fractions seem to demand more working memory resources than representing whole numbers (Halford, 2007). Moreover, fraction knowledge also requires inhibitory control and attention (Siegler et al., 2013), so that the numerator and denominator are not treated as independent whole numbers (Ni, Zhou 2005). With this in mind, it is clear that for a student with MLD (even when "D" stands for "difficulties" instead of "disabilities") the learning of fraction will be a particularly arduous task. In fact, recent studies suggest that dyscalculia, a particular kind of MLD, is rooted specifically in weak visual-spatial working memory and inhibitory control (Szűcs et al., 2013).

Our present work on fractions is part of a broader body of research (Robotti, 2013; Baccaglini-Frank & Robotti, 2013; Baccaglini-Frank, Antonini, Robotti, Santi, 2014) that has the objective of building inclusive curricular material, grounded theoretically in research in mathematics education and in cognitive psychology, appropriate for *all* students, including those with MLD.

**CONCEPTUAL FRAMEWORK**

A large number of studies associated short-term memory (STM) and working memory (WM) with mathematical achievement for students and expert (see reviews in Raghubar, Barnes & Hecht, 2010). Moreover, non-verbal intelligence, addressed to general cognition without reference to the language ability, (DeThorne & Schaefer, 2004), also seems to be strongly related to mathematical achievement (Szűcs et al., 2013). These (and similar) findings suggest that non-verbal intelligence may partially depends on spatial skills (Rourke & Conway, 1997). Thus, spatial processes, performed on the base of spatial skills, can be potentially important in mathematical performances, where explicit or implicit visualization is required. Moreover, research in cognitive science (Stella & Grandi, 2012) has identified specific and preferential channels of access and elaboration of information. For students with MLD these are the visual non-verbal, the kinesthetic and/or the auditory channels.

Studies in mathematics education as well, although with different conceptual frameworks, have highlighted how sensory-motor, perceptive, and kinaesthetic experiences are fundamental for the formation of mathematical concepts – even highly abstract ones (Arzarello, 2006; Gallese & Lakoff, 2005; Nemirovky, 2003; Radford, 2003). In this regard, within a semiotic perspective, Bartolini Bussi and Mariotti (2008) state that the student's use of specific artifacts in solving mathematical problems contributes to his/her development of mathematical meanings, in a potentially "coherent" way with respect to the mathematical meanings aimed at in the teaching activity.

Thus, in this paper we aim to describe examples (activities) of inclusive math education (Ianes & Demo, 2013), constructed referring to the math education domain as well cognitive psychology and neuroscience domains.

The goal of the activities we will describe, was to realize a sequence that would favor, for *all* children (including those with MLD) the development of mathematical meanings of fractions as numbers that can be placed on the number line. The sequence of activities was designed, realized and analyzed taking into account the following principles:

- the importance of an epistemological analysis of the mathematical content
- the role of perceptive and kinesthetic experience in mathematical concept formation as well the visual non-verbal, and auditory channels of access and elaboration of information, in particular in children with MLD
- the role of social interaction, verbalization, mathematical discussion;

- the teacher as a cultural mediator.

Following these principles, as we will describe later, particular artifacts (like paper strips, rulers and scissors) were identified with the intention of using them to help mediate the meanings at stake in the activities.

**METHODOLOGY AND SEQUENCE OF ACTIVITIES**

The sequence of activities was designed by 22 primary school teachers and 1 supervisor (the first author) composing a group of study. The activities were carried out during a pilot experimentation, which involved 22 classes (nine $5^{th}$ grade classes, six $4^{th}$ grade classes and seven $3^{rd}$ grade classes), before being revised for an upcoming full-blown study. In this paper we will report on the pilot experimentation carried out in the $3^{rd}$ grade classes. Students worked in small groups.

The sequence of activities asked to work with differents artifacts such as A4 sheets of paper, squared paper strips or represented squared paper strips in notebooks. At first, were asked to represent fractions on squared paper strips, then it was asked to represent squared paper strips in notebooks and to represent, on them, fractions. At last, students were asked to represent fractions on the number line (see below). As described in the following, the teachers also included moments of institutionalization and discussion based on some significant episodes.

The activities concerning the sequence are:

1) Partitioning of the A4 sheet of paper: this activity involves dividing the A4 sheet of paper, chosen as a unit of measure, in equal parts, by folding and using the ruler; the procedure allows for the introduction of "equivalent fractions" as equivalent surfaces, and of "sum of fractions" for obtaining the whole (the chosen unit, that is, the A4 sheet).
2) Partitioning of a strip of squared paper. This activity involves three sessions:
    A. Given a certain unit of measure, position it on the strip; then, position some fractions on the strip (1/2, 1/3, 1/4, ….) according to the given unit of measure (see Fig. 1). The objective is to represent, on the same strip, different fractions, introducing reciprocal comparison.
    B. Given different units of measure on different strips, on each strip a same fraction is represented (1/2). The objective is to make explicit the dependence of the unit fraction upon the chosen unit of measure (1/2u).
    C. Given a squared strip, choose appropriate units of measure to represent different fractions on that strip (e.g., 1/3 and 1/5). The objective of this activity is to find the lcm (least common multiple) between denominators as the appropriate unit of measure.
3) Placing fractions on the number line. The fractions, considered to be lengths of segments with origin in 0, are placed on the (positive) number line using the idea at the basis of the operator conception of fractions (developed in point 2). Since the right endpoint of the segment on the number line is labeled with fraction, it will also assume the meaning of "number", as all the other whole numbers on the line. Different fractions will be associated to a same point on the line, and will be used to revisit the meaning of "equivalent fractions".

**ANALYSIS AND DISCUSSION**

In this section we present an analysis of points 2 and 3 of the sequence, and in particular the transition from point 2 to point 3, which we consider the most significant in order to place fractions on number line. Our objective is to highlight how the meaning of fraction is evolved, thanks to the use of the tools (squared strip paper, and number line) and to the designed taks.

ACTIVITY 2, SESSION A). A certain unit of measure is given (for instance, a unit measure corresponding to 15 squares). The students are asked to position it on the *strip of squared paper* and to place and color on the strip unit fractions like 1/5,1/3, etc.

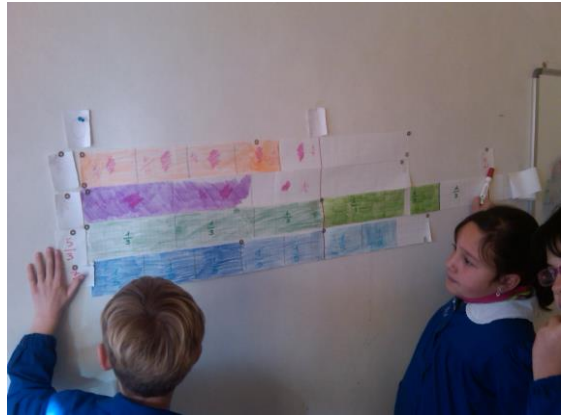

**Fig. 1 Four strips of squared paper where students has postioned a certain unit of measure and had defined unit fractions and colored fractions (4/5, 2/3, 5/3, 7/5)**

With respect to the kinesthetic approach that characterizes activity 1 (partitioning of the A4 sheet of paper in equal part), the manipulation of the artifact "squared strip" becomes a prevalently perceptive experience, in which the main channel for accessing information is the visual non-verbal one. Therefore the task (implicitly) requires the use of a procedure in which the fraction is conceived as an operator: the students partition the strip and produce linguistic signs associated to the name of the fraction expressed in verbal language ("Un mezzo" – tr. "One half"), in verbal visual language (the writing "un mezzo"- tr. One haf) and arithmetical language (1/2). The teacher institutionalizes the relationship between the different signs (partitions of the strips, visual verbal, visual non verbal, and arithmetical signs) in terms of rational numbers. Thus, the construction meaning related to rational number, is based on the interplay among different types of semiotic sets (Arzarello, 2006). Note that the task was completed by all groups of students.

ACTIVITY 2, SESSION B. Each group of Students is asked to choose a unit of measure, reproducing it on a strip and placing the fraction 1/2 on the strip. Then, their strips are compared. The dependence of the fraction on the unit of measure, observed comparing the results of the different groups of students, becomes explicit during a classroom discussion, from which we include an interesting excerpt:

Student 1: Maybe we made a mistake

Student 2: No, we did not make a mistake, I am sure I folded the unit in half, so it's ½

Student 3: We shouldn't look at the length, because each group chose a different unit of measure. […]

Students 4: Because doing ½ is cutting in half, so if the units are different the halves are different […] we have to be careful because to understand which counts more we can't put them one on top of the other like we did for the placemats.

Here a shared meaning is being developed for the fraction as an operator on a chosen unit of measure (1/2u). Note that the kinaestetyc approach (putting strips one biside the other) is no longer an effective strategy in order to compare fractions.

ACTIVITY 2, SESSIONE C). The task asks to choose a unit of measure to represent on the same strip different unit fractions like 1/3, 1/6, 1/8, 1/2, 1/4.

The parameters defining the situations are such that the situation itself makes it evident that it is necessary to choose 24 squares (corresponding to the lcm of 2, 3, 6, 8 and 4) as the unit of measure. In fact the children did not simpy looked for the unit of measure spontaneously, generally using trial and error methods, but they also checked the efficiency of their choice. Moreover, positioning on a single strip different fractions, makes the ordering of fractions exactly like that of the other numbers perceptively evident (Fig 2).

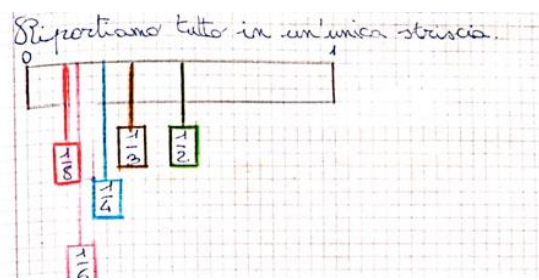 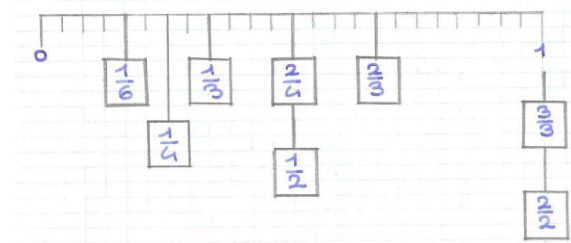

**Fig. 2 different fractions on the same strip    Fig. 3 fractions on the number line**

We note here that for the different fractions on the strips (Fig 2), the teacher asks to associate to the verbal, figural and arithmetical representations also color. The reason is that, as suggested by Stella and Grandi (2011), the verbal channel is not the preferred one for students with MLD. Color becomes a tool supporting working memory and possibly also long term memory, though which the meanings developed can be recalled and used (Baccaglini-Frank, Robotti, 2013).

ACTIVITY 3. The objective of this activity is to place and order fractions on the number line (Fig. 3) that has been partially constructed in Activity 2 and that is now used by the teacher as a tool of semiotic mediation (according to Bartolini Bussi and Mariotti). This transition is fundamental: the representation of the artifact "strip of paper" becomes a mathematical sign that represents the mathematical object "the number line".

Actually, from now on color is no longer used and the labels are referred to points on the number line. We can therefore claim that fractions here have assumed the role of rational numbers. The teacher could take advantage of this transition to construct a new tool of semiotic mediation, developed from the preceding artifacts. The "narrow strip" now becomes a concrete artifact (Fig. 4): it turns into a piece of string on the wall, upon which 0 is placed at the left end and the position of the unit is made to vary dynamically sliding the corresponding label attached with a clothes peg. The dynamic component of this artifact recalls software (such as AlNuSet, GeoGebra, Cabri2…) of course with evident differences, including the fact that as the unit (the position of the paper card with written "1") is made to vary, the positions of the other whole numbers and fractions does not vary dynamically at the same time or automatically, as a consequence of the new placement of the unit: their motion requires a specific action in order to maintain the desired mathematical relationships. This feature can actually be exploited to foster the students' appropriation and *active control* of important mathematical meanings at stake, such as the density of rational numbers or ordering of fractions on numerical line, during activities like this one..

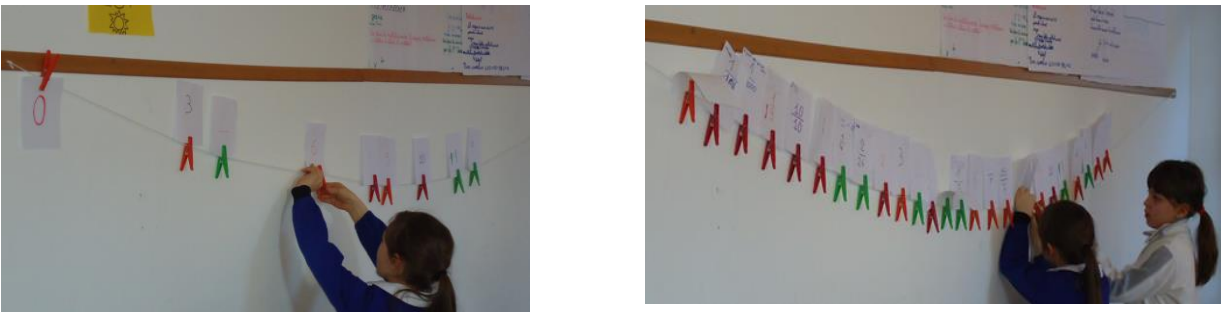

**Fig. 4 String on the wall, the position of the unit is made to vary dynamically sliding the corresponding label attached with a clothes peg**

## CONCLUSIONS

We have outlined particularly significant (and delicate) passages of the sequence of activities, showing how the transition was guided. Initially the students were exposed to a (somewhat traditional) conception of fractions as operators in the context of partitioned areas (the "part-whole" meaning described in Bobis et al., 2013); this idea was soon re-invested in a slightly different context: the areas became strips that gradually lost their "fatness" and were narrowed down until they become (oriented) segments indicating distances from the origin of the number line. The power of an approach like the one described resides in how such a transition can be gradual and continuous, if the teacher manages to keep alive the situated meanings that emerge throughout its unraveling. This is in fact what happened, and the children (including those with certified MLD) came to deal with fractions as numbers on the number line, without hesitating to compare them, place equivalent fractions on the same point, and add them, according to meanings they had developed using the strips of paper that still had an area.

In summary, the analysis of the teaching intervention has showed that students have elaborated personal meanings consistent with the mathematical meanings related to fractions. In particular, the strip was used as instrument of semiotic mediation to develop the meanings related to fractions as operator and, then, to the ordering of fractions, to equivalent fractions and equivalence classes. The use of the strip, the string and the color, has had a key role in favoring the construction of the number line as a mathematical object. On number line, the fractions, associated with points, could assume the role of rational numbers being representatives of equivalence classes. Finally, we think that the construction of meanings related to fractions will also support the management of the procedural aspects involved in operations with fractions, as some researches both in mathematical education domain and cognitive science domain had already advanced (Siegler, 2013 Robotti, 2013, Robotti, Ferrando, 2013). Further studies are needed to explore and to confirm this hypothesis that we consider significant both for research and for teaching.

We would like to say a very great thank you to all teachers of the "Questione di numeri " project who have realized, with one of the author, this research study.